\documentclass[10pt,conference]{IEEEtran}
\usepackage{amsmath}
\usepackage{amsfonts}
\usepackage{amsthm}
\usepackage{graphicx}
\usepackage[capitalize, noabbrev]{cleveref}

\usepackage{algorithm}
\usepackage{algpseudocode}
\usepackage{xcolor}

\newtheorem{theorem}{Theorem}
\newtheorem{lemma}[theorem]{Lemma}

\newcommand{\norm}[1]{\left\| #1 \right\|}

\newcommand{\normsq}[1]{\norm{#1}^2}

\begin{document}

\title{Online Signal Recovery via Heavy Ball Kaczmarz}

\author{   \IEEEauthorblockN{Benjamin Jarman\IEEEauthorrefmark{1}, Yotam Yaniv\IEEEauthorrefmark{1}, Deanna Needell\IEEEauthorrefmark{1}}\IEEEauthorblockA{\IEEEauthorrefmark{1}Department of Mathematics, University of California, Los Angeles \\ \{bjarman, yotamya, deanna\}@math.ucla.edu}}

\maketitle

\begin{abstract}
Recovering a signal $x^\ast \in \mathbb{R}^n$ from a sequence of linear measurements is an important problem in areas such as computerized tomography and compressed sensing. In this work, we consider an online setting in which measurements are sampled one-by-one from some source distribution. We propose solving this problem with a variant of the Kaczmarz method with an additional heavy ball momentum term. A popular technique for solving systems of linear equations, recent work has shown that the Kaczmarz method also enjoys linear convergence when applied to random measurement models, however convergence may be slowed when successive measurements are highly coherent. We demonstrate that the addition of heavy ball momentum may accelerate the convergence of the Kaczmarz method when data is coherent, and provide a theoretical analysis of the method culminating in a linear convergence guarantee for a wide class of source distributions. 
\end{abstract}


\section{Introduction}

\subsection{The Kaczmarz Method}

Recovering a signal $x^\ast \in \mathbb{R}^n$ from a collection of linear measurements is an important problem in computerized tomography \cite{natterer}, sensor networks \cite{savvides2002sensors}, compressive sensing \cite{eldar2012compressed, FoucartCSbook}, machine learning subroutines \cite{Bottou2010LargeScaleML}, and beyond. When the collection of linear measurements is finite, say of size $m$, and accessible at any time, the problem is equivalent to solving a system of linear equations $Ax = b$ with $A \in \mathbb{R}^{m \times n}$ and $b \in \mathbb{R}^m$, which has been well-studied. A popular method for solving this classical problem is the Kaczmarz method \cite{K37:Angena}: beginning with an initial iterate $x_0$, at each iteration a row of the system is sampled and the previous iterate is projected onto the hyperplane defined by the solution space given by that row. More precisely, if the row $a_i^\top x = b_i$ is sampled at iteration $k$, the update has the form
\[
x_{k} = x_{k-1} - \frac{\langle a_i, x_{k-1} \rangle - b_i}{\norm{a_i}^2}a_i.
\]
The original method proposed cycling through rows in order, such that $i = k \text{ mod } m$. In \cite{Herman1993AlgebraicRT} it was observed empirically that randomized row selection accelerates convergence, and in the landmark work \cite{strohmer2009randomized} it was proven that selecting rows at random with probability proportional to their Euclidean norm yields linear convergence in expectation.

In this work, we consider an online model in which at each discrete time $t = 1, 2, \dots$ a linear measurement $(\varphi_t, y_t) \in \mathbb{R}^n \times \mathbb{R}$ is received. We assume that each measurement is noiseless, i.e. $\langle \varphi_t, x^\ast \rangle = y_t$ for all $t$, and that measurements are \emph{streamed} through memory and are not stored. Note that the linear system setting described above is a special case of this model, but we now allow for measurements to be sampled from a more general source. The Kaczmarz method is well-suited to this setting as it requires access to only a single measurement at each iteration. See, for example, \cite{chen2012almost}, where measurement data is viewed as being sampled i.i.d. from some distribution $\mathcal{D}$ on $\mathbb{R}^n$. We assume the noiseless, i.i.d. setting throughout this paper. A Kaczmarz update in this setting has the following form, when initialized with some arbitrary $x_0$: at discrete times $t = 1, 2, \dots$, a measurement $(\varphi_t, y_t) \in \mathbb{R}^n \times \mathbb{R}$ is received, where $\varphi_t \sim \mathcal{D}$, and a Kaczmarz iteration is computed
\[
x_{t} = x_{t-1} - \frac{\langle \varphi_t, x_{t-1}\rangle - y_t}{\norm{\varphi_t}^2}\varphi_t.
\]
In \cite{chen2012almost} it was shown that under certain conditions on $\mathcal{D}$, the method enjoys linear convergence in expectation. Further related works have placed online Kaczmarz in the context of learning theory \cite{lin2015learning}, and have analyzed sparse online variants \cite{lei2018sparse, Lorenz2014ASK}. Random vector models have also appeared in analyses of Kaczmarz methods for phase retrieval \cite{tan2018phase} and for sparsely corrupted data \cite{haddock2020quantilebased}.

\subsection{Heavy Ball Momentum}

Heavy ball momentum is a popular addition to gradient descent methods, in which an additional step is taken in the direction of the previous iteration's movement. Proposed initially in \cite{polyak1964some}, it has proven very popular in machine learning \cite{sutskever2013momentum, krizhevsky2017imagenet,gitman2019understanding,xia2021heavy}, with a guarantee of linear convergence for stochastic gradient methods with heavy ball momentum proven in \cite{loizou2020momentum} (improving on earlier sublinear guarantees in \cite{Yang2016UnifiedCA, gadat2016momentum}).  A gradient descent method itself \cite{needell2016stochastic}, the Kaczmarz method may be modified with heavy ball momentum to give updates of the following form:
\[
x_{t+1} = x_{t} - \frac{\langle \varphi_t, x_{t} \rangle - y_t}{\norm{\varphi_t}^2}\varphi_t + \beta(x_t - x_{t-1}),
\]
where $\beta \geq 0$ is a momentum parameter. In \cite{loizou2020momentum} it was shown that when applied to a linear system (i.e., when each $\varphi_t$ is sampled from the rows of a matrix $A$), the Kaczmarz method with heavy ball momentum converges linearly in expectation. Experimental results indicate accelerated convergence compared to the standard Kaczmarz method on a range of datasets, while the momentum term does not affect the order of the computational cost.

In this work, we propose an \emph{online} variant of the Kaczmarz method with heavy ball momentum. We prove that our method converges linearly in expectation for a wide range of distributions $\mathcal{D}$, and offer particular examples. This theory is supported by numerical experiments on both synthetic and real-world data, which in particular demonstrate the benefits of adding momentum when measurements are highly coherent.

\section{Proposed Method \& Empirical Results}

We propose an online variant of the Kaczmarz method, modified to include a heavy ball momentum term $\beta \in (0,1)$, which we call OHBK($\beta$) (see \cref{alg:hbrk}). We note that our method is a generalization of the momentum Kaczmarz method for systems of linear equations introduced in \cite{loizou2020momentum}. The method requires only a single measurement to be held in storage at a time, while leveraging information about previous measurements through the momentum term.

\begin{algorithm}
\caption{Online Heavy Ball Kaczmarz}\label{alg:hbrk}
\begin{algorithmic}[1]
\Procedure{OHBK($\beta$) }{Input: initial iterate $x_0$, measurements $\{(\varphi_t, y_t)\}_{t=1}^{\infty}$, momentum parameter $\beta$ }
\State Set $x_1 = x_0$
\For{$t = 1, 2, \dots$}
\State Update $x_{t+1} = x_t - \frac{\langle \varphi_t, x_t \rangle - y_t}{\norm{\varphi_t}^2}\varphi_t + \beta(x_t - x_{t-1})$
\EndFor
\EndProcedure
\end{algorithmic}
\end{algorithm}

We test our method on synthetic and real-world data. For each data source, we compare our method OHBK($\beta$) for a variety of $\beta$ to an online Kaczmarz method without momentum, which we denote by OK (equivalently, OHBK($0$)).

\begin{figure}[ht]
    \centering
    \includegraphics[width=\linewidth]{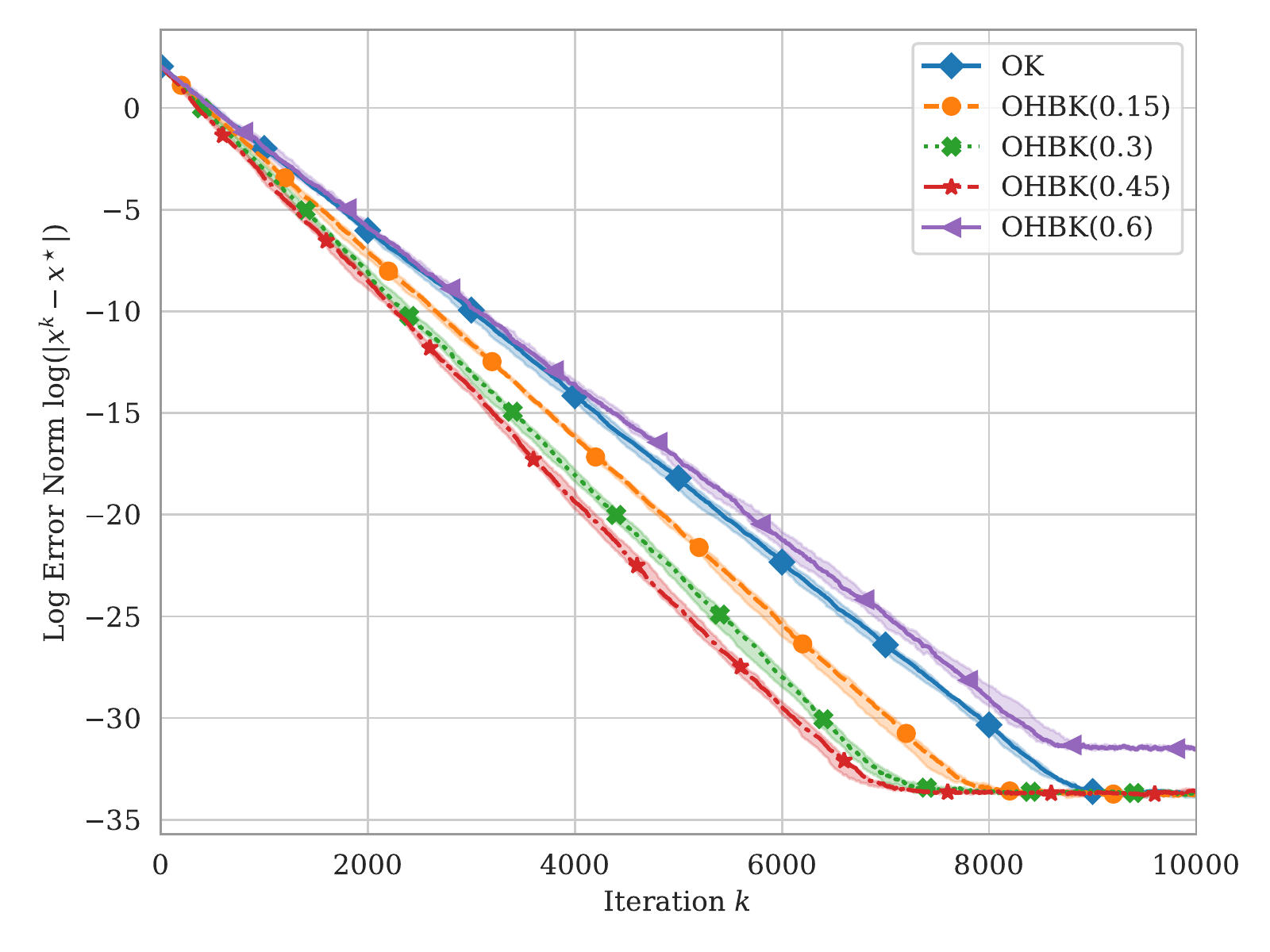}
    \caption{Error versus iteration for OHBK($\beta$) applied to $U[0,1]$ signals of length 50.}
    \label{fig:ohbrk_coherent}
\end{figure}

\begin{figure}[ht]
    \centering
    \includegraphics[width=\linewidth]{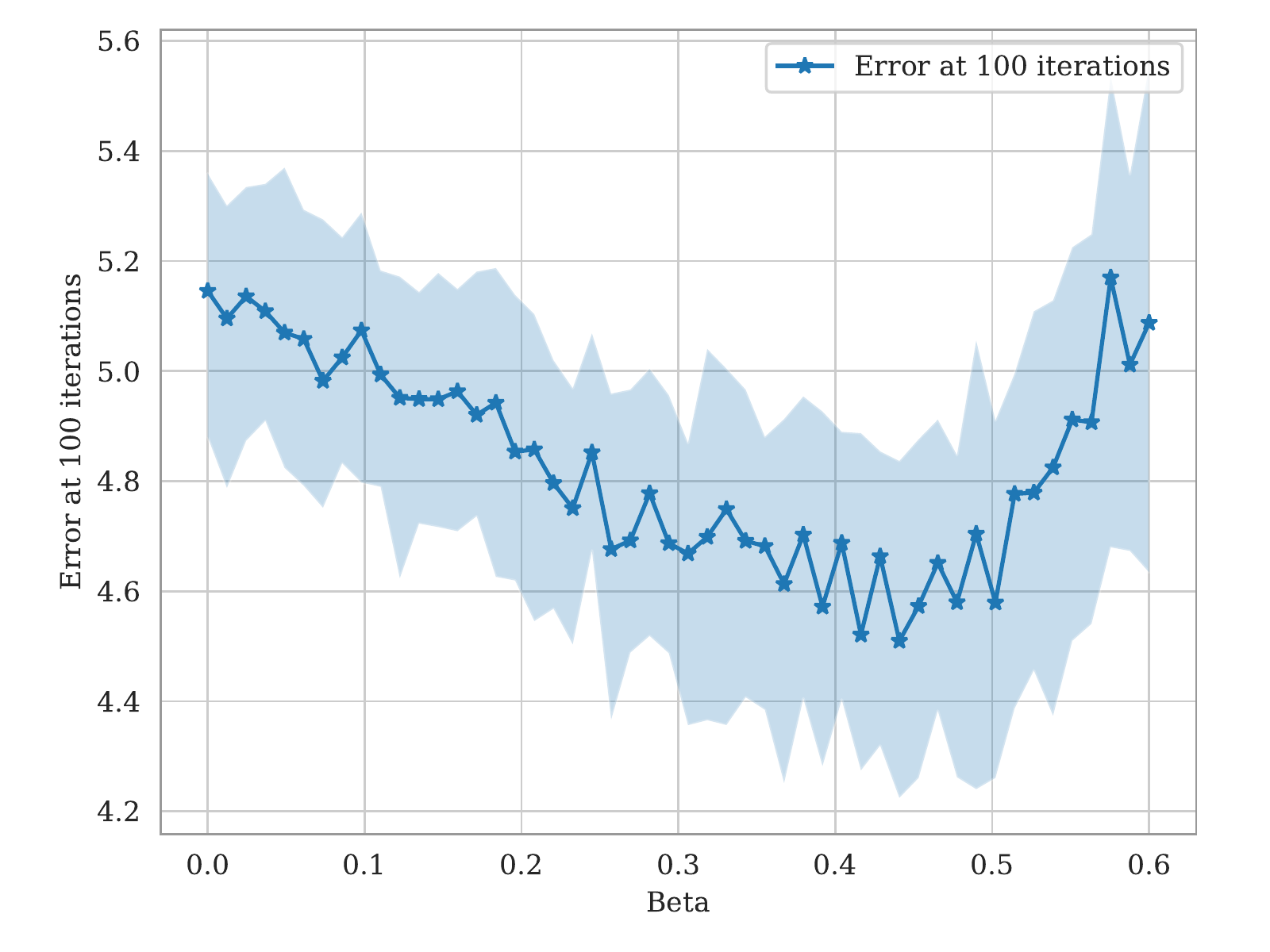}
    \caption{$\norm{x_{100} - x^\ast}$ versus $\beta$ for a range of $\beta \in [0,0.6]$, for $U[0,1]$ signals of length 50.}
    \label{fig:ohbrk_beta_opt}
\end{figure}

We first experiment on synthetic data. We sample $x^\ast \in \mathbb{R}^{50}$ with standard Gaussian entries, and take $\{\varphi_t\}_{t=1}^{\infty}$ to be vectors of length 50 with $U[0,1]$ entries. We note that this process produces particularly coherent data, that is, the vectors $\{\varphi_t\}_{t=1}^{\infty}$ have small pairwise inner products. Each $y_t$ is then computed as $y_t = \langle \varphi_t, x^\ast \rangle$ to ensure measurements are noiseless. In \cref{fig:ohbrk_beta_opt} we perform a parameter search over 100 trials for $\beta$ and plot the median error after 100 iterations versus $\beta$ with shading for the 25th through 75th percentiles. Introducing some amount of momentum provides an acceleration, however, taking $\beta$ to be too large places too much weight on previous information and is less effective. In \cref{fig:ohbrk_coherent} we show convergence down to machine epsilon of OHBK($\beta$) versus online randomized Kaczmarz (i.e. OHBK($0$)) for a selection of $\beta$ (averaging over 10 trials), and the acceleration provided by momentum is clear.

\begin{figure}[ht]
    \centering
    \includegraphics[width=\linewidth]{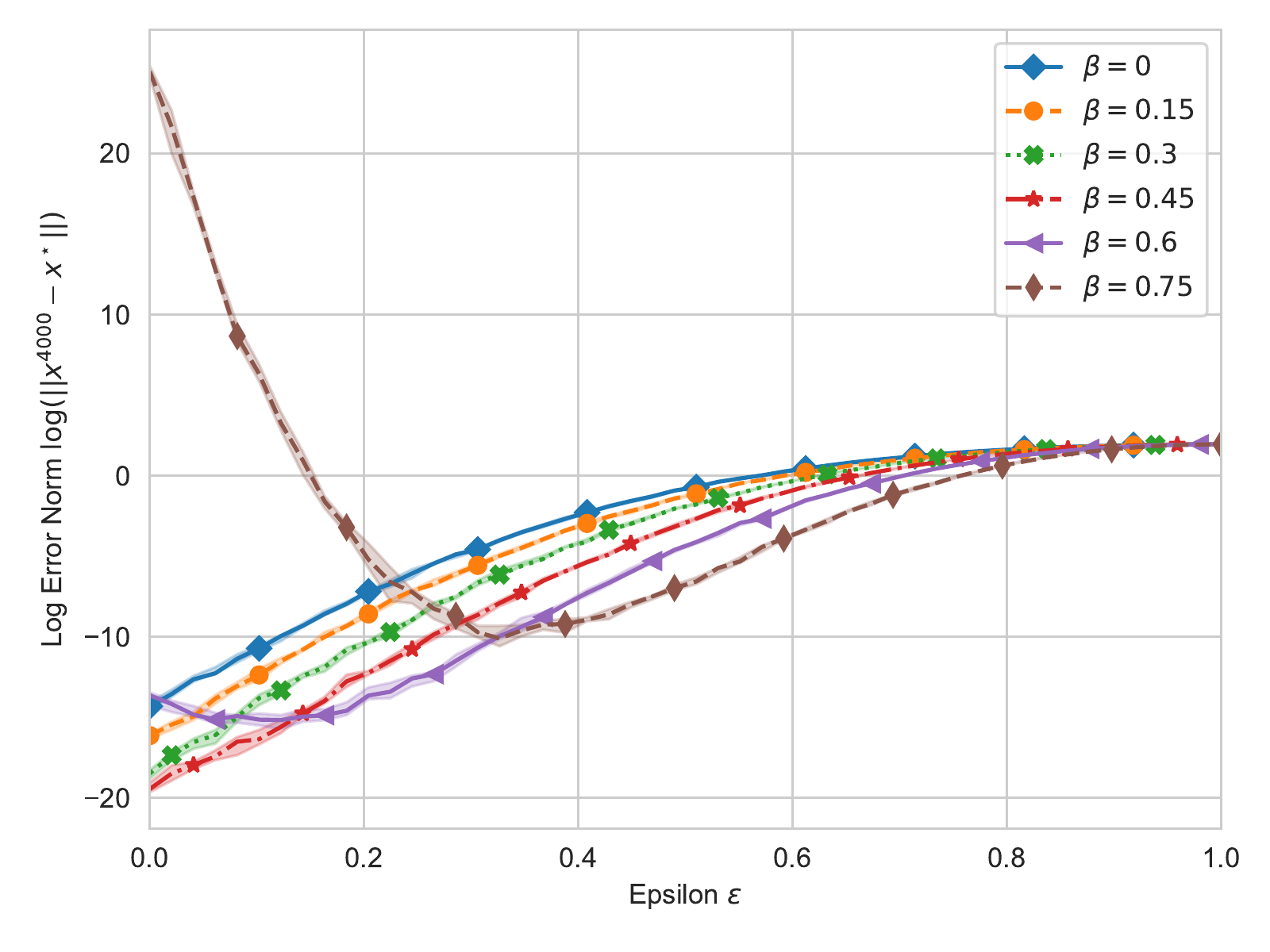}
    \caption{$\log\norm{x_{4000} - x^\ast}$ versus $\varepsilon$ for OHBK($\beta$) applied to $U[\varepsilon, 1]$ signals of length 50.}
    \label{fig:range_of_coherency}
\end{figure}

In \cref{fig:range_of_coherency}, we investigate the effect of momentum on highly coherent systems further. We perform $4000$ iterations of OHBK($\beta$) on $U[\varepsilon, 1]$ signals of length $50$, for $\varepsilon \in [0,1]$, for a range of momentum parameters $\beta$ (again averaged over 10 trials). We see that momentum provides a significant speedup in convergence even for highly coherent systems (i.e. for large $\epsilon$). However, as $\epsilon \to 1$, recovering the signal becomes intractable.

\begin{figure}[ht]
    \centering
    \includegraphics[width=\linewidth]{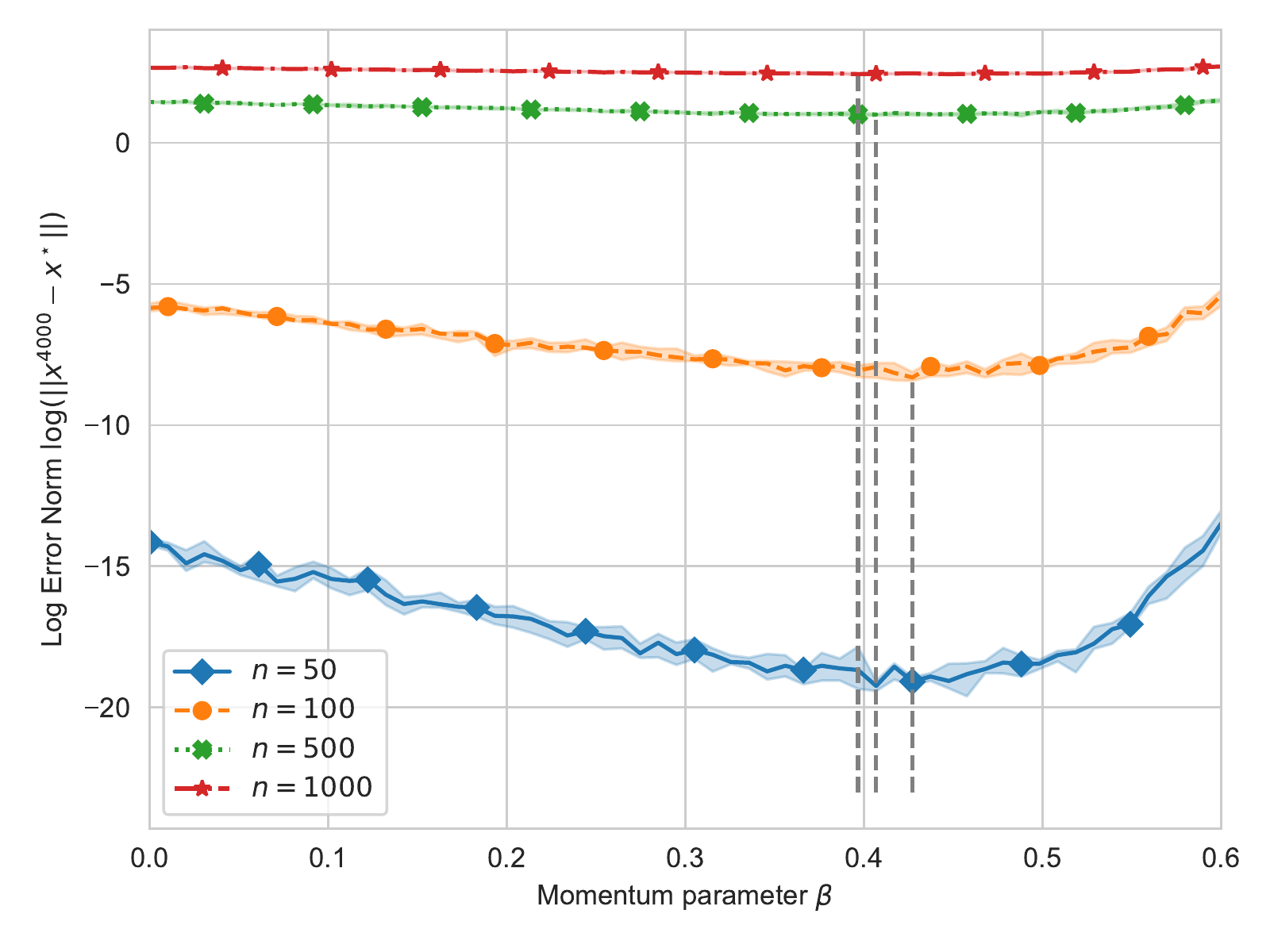}
    \caption{$\log\norm{x_{4000} - x^\ast}$ versus $\beta$ for OHBK($\beta$) applied to $U[0,1]$ signals of length $n$. The gray verticals show the value of $\beta$ yielding the minimum error.}
    \label{fig:range_of_lengths}
\end{figure}

We compare the effect of the signal length $n$ on the optimal momentum parameter $\beta$ in \cref{fig:range_of_lengths}. We perform parameter searches for signals of length $n \in \{50, 100, 500, 1000\}$ and mark the optimal values of $\beta$. The optimal choice of $\beta$ does not appear to vary significantly with $n$.

\begin{figure}[ht]
    \centering
    \includegraphics[width=\linewidth]{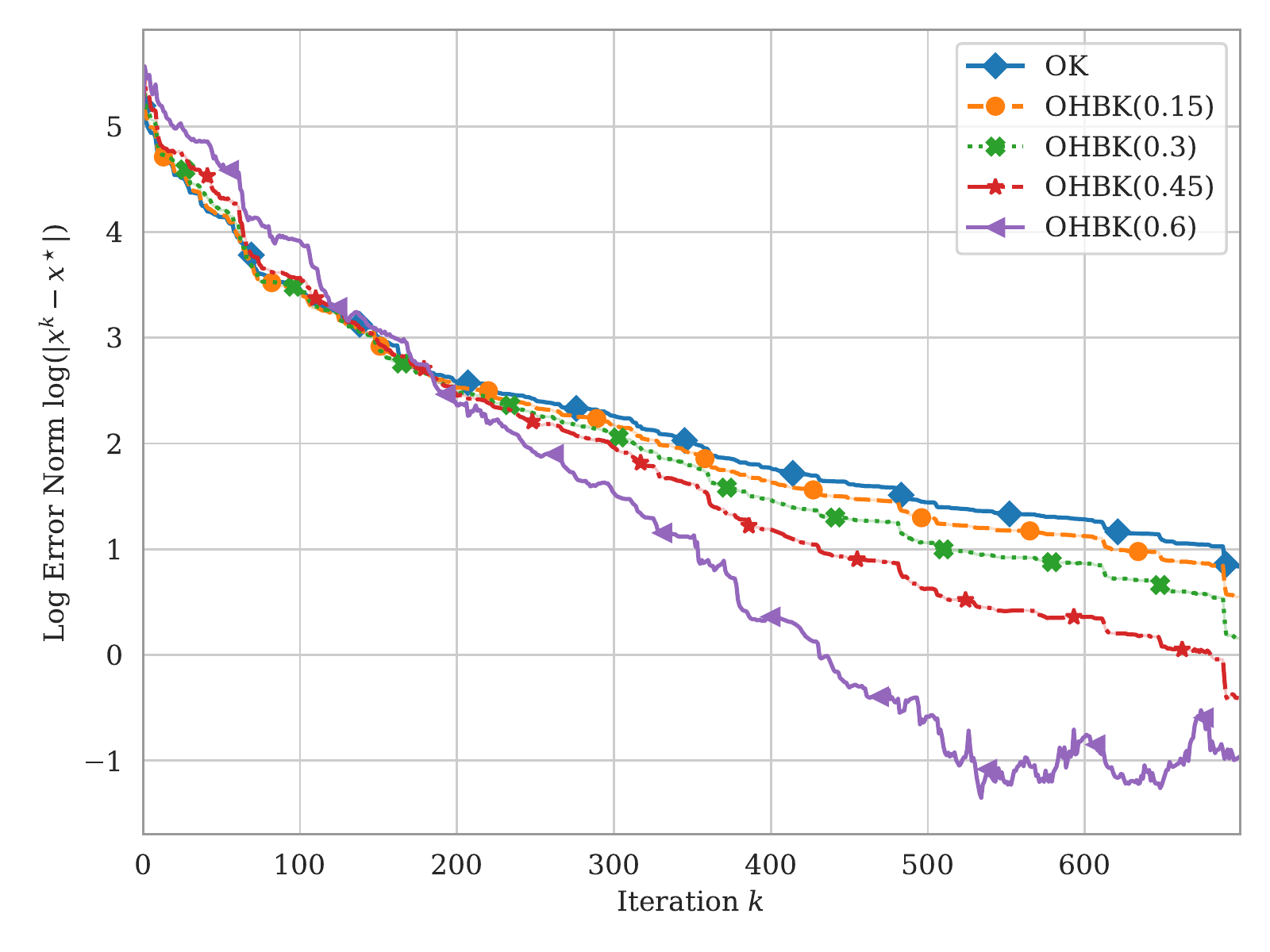}
    \caption{Error versus iteration for OHBK($\beta$) applied to the WDBC dataset.}
    \label{fig:wisco}
\end{figure}

In \cref{fig:wisco} we use a system generated from the Wisconsin Diagnostic Breast Cancer (WDBC) dataset, where each measurement is computed from a digitized image of a fine needle aspirate of a breast mass and describes characteristics of the cell nuclei present \cite{wiscodata}. We stream through each measurement of the 699-row, 10-feature dataset once to replicate the online model, and again see that the addition of momentum provides a noteworthy acceleration to convergence.


\section{Theoretical Results}

Throughout our theory, we assume that $\{\varphi_t\}_{t=1}^{\infty}$ is a sequence of independent samples from some distribution $\mathcal{D}$. We provide a general linear convergence (in expectation) result with a rate depending on the matrix $W := \mathbb{E}_{\mathcal{D}}\left[\frac{\varphi \varphi^\top}{\norm{\varphi}^2}\right]$, in particular on its smallest and largest singular values $\sigma_{\mathrm{min}}(W)$ and $\sigma_{\mathrm{max}}(W)$. 

\begin{theorem}[Convergence in Expectation of OHBRK] \label{thm:1}
Suppose that measurement vectors $\{\varphi_t\}_{t=1}^{\infty}$ are sampled independently from $\mathcal{D}$, and $W = \mathbb{E}_{\mathcal{D}}\left[\frac{\varphi \varphi^\top}{\norm{\varphi}^2}\right]$. Then if $\beta$ is small enough such that
\[
4\beta + 4\beta^2 -(1 + \beta)\sigma_{\mathrm{min}}(W) + \beta \sigma_{\mathrm{max}}(W) < 0,
\]
the iterates produced by OHBK($\beta$) satisfy the following guarantee: for some $\delta > 0$, $q \in (0,1)$, we have
\[
\mathbb{E}[\norm{x_t - x^\ast}^2] \leq q^t (1 + \delta)\norm{x_0 - x^\ast}^2.
\]
\end{theorem}

More interpretable conditions on $\beta$ may be obtained for particular classes of distribution $\mathcal{D}$. In particular, if $\varphi/\norm{\varphi}$ is distributed uniformly on the unit sphere (which is the case if $\mathcal{D}$ itself is the uniform distribution on the unit sphere, or if $\mathcal{D}$ is the standard $n$-dimensional Gaussian), then $W = \frac{1}{n}I$ and we require
\[
\beta + \beta^2 < \frac{1}{4n}
\]
to guarantee linear convergence in expectation. 

\section{Proof of Main Result}
In this section we prove \cref{thm:1} by following the steps of  (\cite{loizou2020momentum}, Theorem 1), making modifications for the online case and simplifications to some of the constants for our special case. First we present a lemma from \cite{loizou2020momentum} which we will use in our convergence proof.

\begin{lemma}[\cite{loizou2020momentum}, Lemma 9]\label{lem:1}
Let $\{F_t\}_{t \geq 0}$ be a sequence of non-negative real numbers with $F_0 = F_1$ that satisfies the relation $F_{t+1} \leq a_1 F_t + a_2 F_{t-1}$ for all $t \geq 1$, with $a_2 > 0$ and $a_1+a_2 < 1$. Then the following inequality hold for all $t \geq 1$ 
\begin{equation*}
    F_{t+1} \leq q^{t}(1+\delta)F_0,
\end{equation*}
where $q = \frac{a_1 + \sqrt{a_1^2+4a_2}}{2} < 1$, $\delta = q-a_1$ and $q \leq a_1+a_2$.
\end{lemma}

A proof of this lemma can be found in \cite{loizou2020momentum}.

We begin our convergence analysis by writing the squared $L2$ error at the $(t+1)$\textsuperscript{th} iteration and substituting the OHBK($\beta$) update into it, 
\begin{align*}
\normsq{x_{t+1} - x^\ast} = \normsq{ x_{t} - \frac{\langle \varphi_t, x_{t} \rangle - y_t}{\norm{\varphi_t}^2}\varphi_t + \beta(x_t - x_{t-1}) - x^\ast}. 
\end{align*}
Next, we group our equation into three terms:
\begin{align}
\begin{split}
    \normsq{x_{t+1} - x^\ast} & = \normsq{x_t - x^\ast - \frac{\langle \varphi_t, x_{t} \rangle - y_t}{\norm{\varphi_t}^2}\varphi_t}\\+&  \beta^2\normsq{x_t - x_{t-1}} \\+& 
    2\beta \langle x_t - x^\ast - \frac{\langle \varphi_t, x_{t} \rangle - y_t}{\norm{\varphi_t}^2}\varphi_t , x_t - x_{t-1}\rangle
\end{split} \label{eq:breakup}
\end{align}

We bound the first term of \cref{eq:breakup} by following a standard Kaczmarz convergence argument and the fact that $y_t = \langle \varphi_t , x^\ast \rangle$. We have that
\begin{align*}
    & \normsq{x_t - x^\ast - \frac{\langle \varphi_t, x_{t} \rangle - y_t}{\norm{\varphi_t}^2}\varphi_t} \\= 
    & \normsq{x_t - x^\ast} + \normsq{\frac{\langle \varphi_t, x_{t} \rangle - y_t}{\norm{\varphi_t}^2}\varphi_t}- \\
    & \qquad 2 \left\langle \frac{\langle \varphi_t, x_{t} \rangle - y_t}{\norm{\varphi_t}^2}\varphi_t, x_t - x^\ast\right\rangle \\ = 
    & \normsq{x_t - x^\ast} + \frac{(\langle \varphi_t, x_{t} \rangle - y_t)^2}{\normsq{\varphi_t}} - 2 \frac{(\langle \varphi_t, x_{t} \rangle - y_t)^2}{\normsq{\varphi_t}} \\=
    & \normsq{x_t - x^\ast} - \frac{(\langle \varphi_t, x_{t} \rangle - y_t)^2}{\normsq{\varphi_t}}.
\end{align*}

We bound the second term of \cref{eq:breakup} by first adding and subtracting $x^\ast$ 
\begin{align*}
   \beta^2\normsq{x_t - x_{t-1}} = \beta^2\normsq{(x_t-x^\ast) + (x^\ast- x_{t-1})}.
\end{align*}
Then by applying the fact that $\normsq{a+b} \leq 2\normsq{a} + 2\normsq{b}$ we have that
\begin{align*}
    &\beta^2\normsq{(x_t-x^\ast) + (x^\ast- x_{t-1})} \\ &\leq 2\beta^2 \normsq{x_t - x^\ast} + 2\beta^2\norm{x_{t-1} - x^\ast}.
\end{align*}

Thus we have that
 \begin{align*}
     \beta^2 \normsq{x_t - x_{t-1}} \leq 2\beta^2 \normsq{x_t - x^\ast} + 2\beta^2\norm{x_{t-1} - x^\ast}.
 \end{align*}

Finally we bound the third term of \cref{eq:breakup} as
\begin{align*}
    & 2\beta \left\langle x_t - x^\ast - \frac{\langle \varphi_t, x_t \rangle - y_t}{\normsq{\varphi_t}}\varphi_t, x_t - x_{t-1}\right\rangle = \\
    & \qquad 2\beta\langle x_t - x^\ast, x_t - x_{t-1}\rangle +\\ & \qquad 2\beta \left\langle \frac{\langle \varphi_t, x_t\rangle - y_t}{\normsq{\varphi_t}}\varphi_t, x_{t-1} - x_t \right\rangle \\
    &= 2\beta \normsq{x_t - x^\ast} + 2\beta \langle x_t - x^\ast, x^\ast - x_{t-1}\rangle + \\
    & \qquad 2\beta \left\langle \frac{\langle \varphi_t, x_t\rangle - y_t}{\normsq{\varphi_t}} \varphi_t, x_{t-1} - x_t\right\rangle \\
    &= \beta \normsq{x_t - x^\ast} + \beta \normsq{x_t - x_{t-1}} - \beta \normsq{x_{t-1} - x^\ast} + \\
    & \qquad 2\beta \left\langle \frac{\langle \varphi_t, x_t\rangle - y_t}{\normsq{\varphi_t}} \varphi_t, x_{t-1} - x_t\right\rangle \\
    &\leq \beta \normsq{x_t - x^\ast} + \beta \normsq{x_t - x_{t-1}} - \beta \normsq{x_{t-1} - x^\ast}- \\
    &\qquad \beta \langle \frac{\langle \varphi_t, x_t\rangle - y_t}{\normsq{\varphi_t}}\varphi_t, x_t - x^\ast\rangle + \\
    &\qquad \beta \langle \frac{\langle \varphi_t, x_{t-1}\rangle - y_t}{\normsq{\varphi_t}}\varphi_t, x_{t-1} - x^\ast\rangle.
\end{align*}

Combining the three bounds, we have 
\begin{align*}
    &\normsq{x_{t+1} - x^\ast} \leq  \normsq{x_t - x^\ast} - \frac{(\langle \varphi_t, x_{t} \rangle - y_t)^2}{\normsq{\varphi_t}} \\+& 2\beta^2 \normsq{x_t - x^\ast} + 2\beta^2\normsq{x_{t-1} - x^\ast} \\ + 
    & \beta \normsq{x_t - x^\ast} + \beta \normsq{x_t - x_{t-1}} - \beta \normsq{x_{t-1} - x^\ast}- \\
    &\qquad \beta \langle \frac{\langle \varphi_t, x_t\rangle - y_t}{\normsq{\varphi_t}}\varphi_t, x_t - x^\ast\rangle +\\
    &\qquad \beta \langle \frac{\langle \varphi_t, x_{t-1}\rangle - y_t}{\normsq{\varphi_t}}\varphi_t, x_{t-1} - x^\ast\rangle.
\end{align*}

Simplifying and grouping like terms we have 

\begin{align*}
    &\normsq{x_{t+1} - x^\ast} \leq  (1+2\beta^2+\beta)\normsq{x_t - x^\ast} +\\ 
    &\qquad (2\beta^2-\beta)\normsq{x_{t-1} - x^\ast} - \\
    & \qquad \frac{(\langle \varphi_t, x_{t} \rangle - y_t)^2}{\normsq{\varphi_t}} +  \beta \normsq{x_t - x_{t-1}} -\\
    &\qquad \beta \langle \frac{\langle \varphi_t, x_t\rangle - y_t}{\normsq{\varphi_t}}\varphi_t, x_t - x^\ast\rangle + \\
    &\qquad \beta \langle \frac{\langle \varphi_t, x_{t-1}\rangle - y_t}{\normsq{\varphi_t}}\varphi_t, x_{t-1} - x^\ast\rangle.
\end{align*}

Applying the simplification for the second term of \cref{eq:breakup} and simplifying the inner products, we have

\begin{align*}
    &\normsq{x_{t+1} - x^\ast} \leq  (1+2\beta^2+3\beta)\normsq{x_t - x^\ast} +\\ 
    &\qquad (2\beta^2+\beta)\normsq{x_{t-1} - x^\ast} - \\
    & \qquad (\beta +1) \frac{\langle \varphi_t, x_{t} -x^\ast \rangle^2}{\normsq{\varphi_t}} +\\
    & \qquad \beta \frac{\langle \varphi_{t}, x_{t-1} - x^\ast \rangle^2}{\normsq{\varphi_{t}}}.
\end{align*}

Taking an expectation over our signal of our simplified equation

\begin{align*}
    &\mathbb{E}[\normsq{x_{t+1} - x^\ast}] \leq  (1+2\beta^2+3\beta)\normsq{x_t - x^\ast} +\\ 
    &\qquad (2\beta^2+\beta)\normsq{x_{t-1} - x^\ast} - \\
    & \qquad (\beta +1) \mathbb{E}[\frac{(\langle \varphi_t, x_{t} -x^\ast \rangle)^2}{\normsq{\varphi_t}}] +\\
    & \qquad \beta \mathbb{E}[\frac{(\langle \varphi_{t}, x_{t-1} - x^\ast \rangle)^2}{\normsq{\varphi_{t}}}] \\
    &= (1+2\beta^2+3\beta)\normsq{x_t - x^\ast} +\\ 
    &\qquad (2\beta^2+\beta)\normsq{x_{t-1} - x^\ast} - \\
    & \qquad (1+\beta)(x_t - x^\ast)^T \mathbb{E}\left[\frac{\varphi_t \varphi_t^T}{\normsq{\varphi_t}}\right](x_t - x^\ast) + \\
    & \qquad \beta(x_{t-1} - x^\ast)^T \mathbb{E}\left[\frac{\varphi_t \varphi_t^T}{\normsq{\varphi_t}}\right](x_{t-1} - x^\ast).
\end{align*}

Let $W := \mathbb{E}\left[\frac{\varphi_t \varphi_t^T}{\normsq{\varphi_t}}\right]$. We can then bound the above in terms of the largest and smallest singular values of $W$:

\begin{align*}
    &\mathbb{E}[\normsq{x_{t+1} - x^\ast}] \leq  (1+2\beta^2+3\beta)\normsq{x_t - x^\ast} +\\ 
    &\qquad (2\beta^2+\beta)\normsq{x_{t-1} - x^\ast} - \\
    & \qquad (1+\beta)\sigma_{\min}(W)\normsq{x_t - x^\ast} + \\
    & \qquad \beta\sigma_{\max}(W)\normsq{x_{t-1} - x^\ast} \\
    &= (1+2\beta^2+3\beta - (1+\beta)\sigma_{\min}(W))\normsq{x_t - x^\ast} + \\
    & \qquad (2\beta^2+\beta + \beta\sigma_{\max}(W))\normsq{x_{t-1} - x^\ast}.
\end{align*}

Finally, we apply \cref{lem:1}, wherein the two coefficients are given by $a_1 = 1+2\beta^2+3\beta - (1+\beta)\sigma_{\min}(W)$ and $a_2 = 2\beta^2+\beta + \beta\sigma_{\max}(W)$. Since we assumed that $a_1 +a_2 = 1+4\beta^2+4\beta +(1+\beta)\sigma_{\min}(W) + \beta \sigma_{\max}(W) < 1$ and since $\beta > 0$ then $a_2 = 2\beta^2+\beta + \beta\sigma_{\max}(W) > 0$ thus the assumptions for \cref{lem:1} hold, so we have that 
\[
\mathbb{E}[\norm{x_t - x^\ast}^2] \leq q^t (1 + \delta)\norm{x_0 - x^\ast}^2
\]

where $q = \frac{a_1 + \sqrt{a_1^2+4a_1}}{2}$, $\delta = q - a_1$ and $a_1+a_2 \leq q < 1$. Since $q \in (0,1)$ we have shown that the norm squared error of the iterates produced by OHBK($\beta$) converges linearly in expectation. 

\section{Conclusion and Future Directions}
In this work we discuss using a Kaczmarz method variant with momentum to solve an online signal recovery problem. We leverage a heavy ball momentum term, a classical acceleration method, to improve the convergence rate. We prove a theoretical convergence rate for OHBK($\beta$), and verify this convergence empirically on both synthetic and real-world data. We demonstrate empirically that for coherent measurements, the addition of momentum indeed accelerates convergence, and provided some initial exploration into the dependence of the convergence rate on the signal length $n$ and momentum strength $\beta$. 

It is notable that in our convergence analysis, we did not recover a theoretically optimal value for $\beta$. Doing so, and comparing this value to empirically best values, would be an interesting future direction. Furthermore, we would like to obtain theoretical parameter relationships: for example, how the optimal momentum strength depends on the signal length and coherency of the measurements. It may in fact be optimal to adaptively adjust the momentum parameter across iterations based on the current iterate and properties of incoming measurements. Additionally, we would like to leverage other accelerated gradient methods such as ADAM \cite{kingma2014adam}. Finally, we would like to consider solving the online signal recovery problem in the case where each measurement is no longer exact, but instead contains some amount of noise \cite{needell2010randomized}. This could be achieved, for example, using relaxation.

\section*{Acknowledgments}

BJ and DN were partially supported by NSF DMS-2108479, YY and DN were partially supported by NSF DMS-2011140.

\bibliographystyle{IEEEtran}
\bibliography{IEEEabrv,ref}

\end{document}